\documentclass{amsart}

\newtheorem{theorem}{Theorem}[section]
\newtheorem{lemma}[theorem]{Lemma}

\theoremstyle{definition}

\newtheorem{example}[theorem]{Example}

\newcommand{\bF}{\bar{F}}
\newcommand{\bL}{\bar{L}}

\newcommand{\bV}{\bar{V}}

\newcommand{\C}{\mathcal C} 
\newcommand{\D}{\mathcal D} 
\newcommand{\p}{{\mbox{$[p]$}}}
\newcommand{\F}{\mathbb{F}}
\newcommand{\scp}{{\mbox{$\scriptstyle [p]$}}}

\DeclareMathOperator{\cl}{Cl}
\DeclareMathOperator{\ind}{Ind}
\DeclareMathOperator{\Hom}{Hom}
\DeclareMathOperator{\End}{End}
\DeclareMathOperator{\res}{Res}

\numberwithin{equation}{section}

\title[Induced modules]{Induced modules for modular Lie algebras }
\author{Donald W. Barnes}
\address{1 Little Wonga Rd.\\Cremorne NSW 2090\\Australia\\}
\email{donaldbarnes@tpg.com.au}

\thanks{It has been pointed out to me that the results of this paper all appear in Farnsteiner \cite{Farn}. His priority is acknowledged. }

\subjclass[2010]{Primary 17B50}
\keywords{Lie algebras,   induced modules}
\begin{document}

\begin{abstract} Let $L$ be a finite-dimensional Lie algebra over a field of non-zero characteristic and let $S$ be a subalgebra.  Suppose that $X$ is a finite set of finite-dimensional $L$-modules. Let $\D$ be the category of all finite-dimensional $S$-modules. Then there exists a category $\C$ of finite-dimensional $L$-modules containing the modules in $X$ such that the restriction functor $\res_\D^\C: \C \to \D$ has a left adjoint $\ind_\D^\C: \D \to \C$.
\end{abstract}
\maketitle
\section{Introduction} \label{intro}
In the theory of finite groups, much use is made of induced representations.  For an account of this theory and some examples of its use, see Curtis and Reiner \cite[Chapter VIII]{CR}.
If $G$ is a finite group, $H$  a subgroup and $W$ is an $FH$-module, then the induced module $\ind_H^G(W)$ is defined to be $FG \otimes_{FH}W$.  Here, $FG$ denotes the group algebra of $G$ over the field $F$.  An important property of this construction is Frobenius reciprocity.  This, in modern terminology, essentially is that induction is a left adjoint (see Mac Lane \cite{McL} or the anonymous article \cite{Wiki}) to the restriction functor $\res_H^G$, the functor which converts $FG$-modules into $FH$-modules by forgetting the action of elements of $G$ which are not in $H$, that is, that there exists a natural isomorphism $\Hom(\ind_H^G( W), V) \to \Hom(W, \res_H^G(V))$.  If in this isomorphism, we take $V = \ind_H^G(W)$, then corresponding to the identity $\ind_H^G(W) \to \ind_H^G(W)$, we get a monomorphism 
$W \to \res_H^G(\ind_H^G(W))$, known as the unit of the adjunction.

For a Lie algebra $L$, subalgebra $S$ and $S$-module $W$, using the universal enveloping algebras $U(L),  U(S)$,  we can construct an induced $L$-module $U(L) \otimes_{U(S)} W$.  This construction gives a left adjoint to the restriction functor, but is of little use in the theory of finite dimensional Lie algebras since, apart from trivial cases, the induced module so constructed is infinite dimensional.  There have been several attempts to define induced modules which are finite dimensional.  Of necessity, these involve either restrictions on the modules covered by the theory or a weakening of the requirement that the induction functor be a left adjoint to the restriction functor.

For Lie algebras $L$ over a field $F$ of characteristic $0$, in Zassenhaus \cite{Z} and in Hochschild and Mostow \cite{HM}, a construction is given for the special case where $S$ is an ideal of $L$ and $L$ is the vector space direct sum of $S$ and a subalgebra $H$ of $L$.  This construction is not a left adjoint to the restriction functor, but it does have the property of the unit of an adjunction mentioned above.   The module constructed from $W$ contains an $S$-submodule isomorphic to $W$.  This construction is used in \cite{HM} to prove a strengthening of Ado's Theorem.

In \cite{W1, W2}, Wallach gives a construction essentially limited to the case where $F$ is algebraically closed of characteristic $0$, $L$ is semisimple and $S$ is a Borel subalgebra.  Wallach's construction is based on the coinduced construction, $\Hom_{U(S)}(UL),W)$ rather than the tensor product.  Again, the construction has a weakened form of some of the properties of a left adjoint.  Using it, Wallach gave a construction for all irreducible modules of a semisimple Lie algebra.

Over a field $F$ of characteristic $p$, properties of restricted Lie algebras and of $p$-envelopes open  up possibilities not available in characteristic $0$.
For a restricted Lie algebra $(L,\p)$, $\p$-subalgebra $S$ and character $c:L \to F$, we can use the $c$-reduced enveloping algebra $u(L,c)$ to define, for an $S$ module $W$ with character $c|S$, the $c$-induced module $u(L,c) \otimes_{u(S,c|S)}W$.  See Strade and Farnsteiner \cite[Sections 5.3, 5.6]{SF}.  This does give a left adjoint to $\res: \C \to \D$ where $\C$ is the category of $L$-modules with character $c$ and $\D$ the category of $S$-modules with character $c|S$.  Not every $S$-module has a character.  Only in the case where the field $F$ is algebraically closed and the module $W$ is irreducible, is the existence of a character guaranteed.  Despite this limitation, this construction, combined with the use of $p$-envelopes, is used to prove that, for a soluble Lie algebra over an algebraically closed field of characteristic $p$, every irreducible module has dimension a power of $p$.

In Barnes \cite{Cluster}, the concept of a simple character cluster was introduced and used to define  cluster-induced modules.  Three conditions were required for this construction to work:  the cluster $C$ had to restrict simply to $S$, that is, distinct characters in $C$ had to have distinct restrictions to $S$, the field $F$ had to be perfect and the $S$-module $W$ had to be amenable, that is, the module $\bar{V}$ obtained by extending the field to the algebraic closure had to be a direct sum of modules with character in $C|S$.  In this paper,  constructions are given which do not require these conditions.  

If $\C$ is a category of $L$-modules, we have the restriction functor $\res: \C \to \D$ where $\D$ is a category of $S$-modules containing at least those obtained from modules in $\C$ by restriction of the action.  In this paper, we construct an induction functor left adjoint to $\res$ for some categories $\C, \D$ of finite dimensional modules.    We shall see that for any given finite set $X$ of $L$-modules and set $Y$ of $S$-modules, we can choose $\C, \D$ containing the given modules.  The dimension of the induced module $\ind_\D^\C(W)$ constructed depends on the choice of the category $\C$ but not on the choice of $\D$.  We may always choose $\D$ to be the category of {\em all} finite-dimensional $S$-modules.

In Section \ref{restricted}, we give the construction for induced modules for restricted Lie algebras.  In Section \ref{examples}, we give some illustrative examples and see some ways in which  in some cases, the construction may be modified to reduce the dimension of the constructed modules.  In Section \ref{envelopes}, we extend the construction to non-restrictable Lie algebras.

\section{Generalised induction for restricted algebras.} \label{restricted}

In the following, $(L,\p)$ is a restricted Lie algebra over the field $F$, $\bF$ is the algebraic closure of $F$ and $\bL = \bF \otimes_F L$ is the algebra obtained by extension of the field.  As we will have $L$ embedded in various associative algebras, we will denote the Lie algebra product of $x,y$ by $[x,y]$, reserving the notation $xy$ for the product in the associative algebra under consideration.  

 By a character of $L$ we shall understand an $F$-linear map $c : L \to \bF$.  If $\{e_1, \dots, e_n\}$ is a basis of $L$, then $c$ can be expressed as a linear form $c(x) = \sum a_ix_i$ for $x = \sum x_ie_i$, where $a_i \in \bF$.  If $\alpha$ is an automorphism of $\bF/F$, that is, an automorphism of $\bF$ which fixes all elements of $F$, then $c^\alpha$ is the character $c^\alpha(x) = \sum a_i^\alpha x_i$ and is called a conjugate of $c$.   We do not distinguish in notation between $c:L \to \bF$ and its linear extension $\bL \to \bF$.  A character cluster is a set $C$ of characters  which contains, along with any character $c$, all conjugates $c^\alpha$ of $c$.  

If $V$ is an $L$-module, then $\bV$ is the $\bL$-module $\bF \otimes_F V$.  The action of $x \in L$ on $V$ is denoted by $\rho(x)$.  The module $V$ has character $c$ if $(\rho(x)^p-\rho(x^\scp))v = c(x)^pv$ for all $x \in L$ and all $v \in V$.  The cluster $\cl(V)$ of an $L$-module $V$ is the set of all characters of composition factors of $\bar{V}$.

In the universal enveloping algebra $U(L)$, the element $z_x=x^p -x^\scp$ is central.  (See Strade and Farnsteiner \cite[p.203]{SF}.)  For the module $V$ giving the representation $\rho$, we put $\phi_x = \rho(x)^p - \rho(x^\scp)$.  We then have $[\phi_x, \rho(y)] = 0$ for all $x,y \in L$.  By \cite[Lemma 2.1]{Cluster}, 
the map $\phi : L \to \End(V)$ defined by $\phi_x(v) = (\rho(x^p) - \rho(x^\scp))v$ is $p$-semilinear.

Let $W$ be an $S$-module.  To use the standard tensor product construction of induced modules, we need enveloping algebras $u(L), u(S)$ which are quotients of the universal enveloping algebras $U(L),U(S)$ such that 
\begin{enumerate} 
\item $W$ is a $u(S)$-module, 
\item $u(S)$ is a subalgebra of $u(L)$, and
\item $u(S)$ has finite codimension in $u(L)$.
\end{enumerate}
Condition (1) is easily satisfied.  As $S$ is a $U(S)$-module, we can take $K$ the kernel of the representation of $U(S)$ on $W$ and put $u(S) = U(S)/K$.  It is in order to make choices which also satisfy conditions (2) and (3) that we use the extra structure of restricted Lie algebras.  We modify the construction of character reduced enveloping algebras given in Strade and Farnsteiner \cite[Section 5.3]{SF}.

Let $\{e_i \mid i \in I\}$ where $I$ is a finite ordered index set, be an ordered basis of $L$.  We use the multi-index notation of \cite[p. 51]{SF}.  A product of $0$ factors is interpreted as the element $1 \in U(L)$.  The elements $z_i = e_i^p - e_i^\scp$ are in the centre of $U(L)$.  We use the standard filtration of $U(L)$,  $U_k = \langle x_1 \dots x_r \mid  x_j \in L, r\le k \rangle$.  Let $f$ be a family of polynomials $f_i(t) \in F[t]$  of degrees $d_i \ge 1$.  Putting $\zeta_i = f_i(z_i)$ and $v_i = e_i^{pd_i} - \zeta_i$, we have that $\zeta_i$ is in the centre of $U(L)$ and that $v_i \in U_{pd_i - 1}$.  By \cite[Lemma 1.9.7]{SF}, $\{e^\alpha \zeta^\beta \mid \alpha_i < pd_i\}$ is a basis of $U(L)$.  Let $K$ be the ideal of $U(L)$ generated by the $e^\alpha \zeta^\beta$ with $\beta \ne 0$ and let $u(L, f) = U(L)/K$.  It follows that $\{e^\alpha + K \mid  \alpha_i < p d_i\}$ is a basis of $u(L,f)$.  We call $u(L,f)$ the $f$-reduced enveloping algebra of $L$.  It is clearly finite-dimensional.

We now choose the basis and the $f_i$ to construct  algebras $u(L), u(S)$ satisfying the conditions (1), (2) and (3) above.  We partition $I = I_1 \cup I_2$ with the elements of $I_2$ coming after those of $I_1$ and choose a basis such that $\{e_i \mid i \in I_2\}$ is a basis of $S$.  We choose the polynomials $f_i$ for $i \in I_2$ such that $f_i(z_i)W = 0$.  It follows immediately that $\{e^\alpha + (K \cap U(S)) \mid \alpha_i = 0 \text{ for } i \in I_1, \alpha_i < pd_i \text{ for } i \in I_2\}$ is a basis of $u(S) = U(S)/K\cap U(S)$, the $f|I_2$-reduced enveloping algebra of $S$.  By our choice of the $f_i$ for $i \in I_2$, $K \cap U(S)$ is in the kernel of the representation of $U(S)$ on $W$, so $W$ is a $u(S)$-module.  Further, $u(S)$ is a subalgebra of $u(L) = u(L,f)$.  We call these algebras $u(L), u(S)$ the $f$-reduced enveloping algebras and denote them by $u(L,f)$ and $u(S,f|S)$.  Using them, we construct the $f$-induced module $\ind_S^L(W,f) = u(L,f) \otimes_{u(S,f|S)} W$.  Let $\C$ be the category of all finite-dimensional $u(L,f)$-modules and $\D$ be the category of all finite-dimensional $u(S,f|S)$-modules.  Then $\ind_S^L(\_,f): \D \to \C$ is a left adjoint to the restriction functor $\res^L_S: \C \to \D$.

The dimension of the induced module is $\dim(W)\Pi_{i \in I_1}(pd_i)$.  Note that the induced module $\ind_S^L(W,f)$ is not affected by the choice of the $f_i$ for $i \in I_2$, subject to the requirement that $f_i(z_i)W = 0$.  We do not need to choose $f_i$ to be the minimum polynomial.  If several modules are under consideration, we can choose any $f_i$ such that $f_i(z_i)$ annihilates all those modules.  This leads us to modify the above construction by allowing the polynomials $f_i$ for $i \in I_2$ to be zero.  If any $f_i=0$, then $u(L,f)$ and $u(S,f|S)$ are infinite-dimensional, but the codimension remains finite and the induced module is unaffected by this change.  If we take $f_i=0$ for all $i \in I_2$, then $u(S,f|S) = U(S)$ and $\D$ is the category of all finite-dimensional $S$-modules.

\begin{theorem} \label{Thres} Let $(L,\p)$ be a restricted Lie algebra and let $S$ be a \p-subalgebra of $L$.  Let $X$ be a finite set of finite-dimensional $L$-modules.  Let $D$ be the category of all finite-dimensional $S$-modules.  Then there exists a category $\C$ of finite-dimensional $L$-modules such that $X\subseteq  \C$ and the restriction functor $\res:\C \to \D$ has a left adjoint $\ind_\D^\C$.
\end{theorem}

\begin{proof}  Let $V$ be the direct sum of the modules in $X$ and let $\rho$ be the representation of $L$ on $V$.  We choose an ordered basis $\{e_i \mid i \in I\}$ of $L$, partitioned $I = I_1 \cup I_2$ as above.  For $i \in I_1$, we take for $f_i$ the minimal polynomial of $\rho(e_i^p - e_i^\scp)$.  For $i \in I_2$, we take $f_i=0$.  Let $\C$ be the category of all finite-dimensional $u(L,f)$-modules.  Then $\ind_\D^\C = \ind_S^L(\_,f) : \D \to \C$ is the required left adjoint.
\end{proof}

\section{Illustrative examples} \label{examples}

We begin by applying the above construction to the module of Example 6.6 of \cite{Cluster} and comparing the result with the module constructed in that paper.  
\begin{example} \label{ex1}
Let $F$ be the field of three elements, $L = \langle x,y\mid [x,y]=y \rangle, x^\scp = x, y^\scp = 0$,  $S = \langle x \rangle$  and $W = \langle b^1, b^2 \rangle$ with $xb^1 = b^2$ and $xb^2 = -b^1$.  The cluster  used was $C = \{c_1,c_2 \}$ with $c_1(x) = i, c_2(x) = -i$ and $c_1(y) = \alpha + \beta i, c_2(y) = \alpha - \beta i$ where $\alpha, \beta \in F$ and $i^2 = -1$.  Consider the case where $\beta \ne 0$.

 To apply the construction of Section \ref{restricted}, we take the ordered basis $e_1 = y, e_2=x$  and use the polynomial $f_2(t) = 0$.  In order that the induced module shall have the characters $c_1,c_2$ above in its cluster, we take $f_1(t) = t^2 +\alpha t + \alpha^2 + \beta^2$ to produce the following module which we compare to Example 6.6 of \cite{Cluster}.

Consider the module $M = \ind_S^L(W,f)$.  This has basis the twelve elements $m_j^r = y^r \otimes b^j$ for $r = 0,1, \dots, 5$ and $j = 1,2$.   
In the algebra $u(L,f)$, $f_1(y^3-y^\scp)=0$.  Since $y^\scp = 0$, this gives $y^6 +\alpha y^3 +\alpha^2 +\beta^2 = 0$.  Also in $u(L,f)$, we have $xy = y + yx, xy^2= -y^2 + y^2x$, $xy^3 = y^3x$,
$xy^4 =  y^4 +y^4x$, $xy^5 =  -y^5+y^5x$. Using these equations we can now calculate the action of $L$ on $M$.  

\begin{align*}
xm_1^0 &= x(1 \otimes b^1) = m_2^0, & xm_2^0&=x(1\otimes b^2)=-m_1^0,\\  
xm_1^1 &= (y+yx)\otimes b^1 =m_1^1+m_2^1, & xm_2^1&=(y+yx)\otimes b^2=m_2^1 - m_1^1,\\ 
xm_1^2&=(-y^2+y^2x)\otimes b^1= -m_1^2 + m_2^2, & xm_2^2&= (-y^2+y^2x)\otimes b^2=-m_2^2-m_1^2,\\
\end{align*}
\begin{align*}
xm_1^3&=y^3x \otimes b^1= m_2^3,&xm_2^3&=-m_1^3,\\
xm_1^4&= (y^4+y^4x)\otimes b^1 = m_1^4 +m_2^4,&xm_2^4&= -m_1^4 + m_2^4, \\
xm_1^5&= (-y^5+y^5x) \otimes b^1 = -m_1^5+m_2^5, &xm_2^5&=-m_1^5  -m_2^5,\\
ym_1^0&=m_1^1, & ym_2^0&=m_2^1,\\
ym_1^1&=m_1^2, & ym_2^1&=m_2^2,\\
ym_1^2&=m_1^3, & ym_2^2&=m_2^3,\\
ym_1^3&=m_1^4, &ym_2^3&=m_2^4,\\
ym_1^4&=m_1^5, &ym_2^4&=m_2^5,\\
ym_1^5&=-(\alpha y^3 +\alpha^2 +\beta^2) \otimes b^1 
&ym_2^5&=-(\alpha y^3 +\alpha^2 +\beta^2) \otimes b^2 \\
&=-(\alpha^2+\beta^2)m_1^0 -\alpha m_1^3,&&=-(\alpha^2+\beta^2)m_2^0 -\alpha m_2^3.
\end{align*}
\end{example}
Unlike the calculation in \cite[Example 6.6]{Cluster} of the action on $\ind_S^L(W,C)$, the calculation of the action of $L$ on $\ind_S^L(W,f)$ does not require use of any extension of the field $F$.  But $\dim(\ind_S^L(W,f)) = 12$, while $\dim(\ind_S^L(W,C))=6$.  The reason for this is that there are two ways the selected character values for $y$ may be combined with those for $x$, giving two different clusters $C_1, C_2$ with $\ind_S^L(W,f) = \ind_S^L(W,C_1) \oplus \ind_S^L(W,C_2)$.  We can recover $\ind_S^L(W,C_1)$ from $\ind_S^L(W,f)$ by making $\C$ the category of finite-dimensional $u(L,f)$-modules with cluster $C_1$ and taking for  $\ind_\D^\C(W)$  the $C_1$-component of $\ind_S^L(W,f)$.

In the case $\beta=0$, we can take $f_1(t)=t-\alpha$ and we obtain $\ind_S^L(W,f)$ isomorphic to the module  $\ind_S^L(W,C)$ of \cite[Example 6.6]{Cluster} for this case.  If, however, we use $f_1(t) = (t-\alpha)^2$, then $\ind_S^L(W,f)$ is indecomposable with two composition factors, both isomorphic to $\ind_S^L(W,C)$.

Finally, we give an example illustrating the construction when the field is not perfect.

\begin{example} Let $F = \F_3(\tau)$ be the field of rational functions in the indeterminate $\tau$ over the field $\F_3$  of $3$ elements.  Again we take $L = \langle x,y \rangle$ with $[x,y]=y$ and $x^\scp = x$, $S = \langle y \rangle$ and $W = \langle w \rangle$ with $yw = w$.  We have to choose a character $c: L \to F$.  Only the value of $c(x)$ affects the induced module constructed. The minimal polynomial of $x^3-x^\scp$ over $\bF$ is $t-c(x)^3$, so to get a different result from the perfect field case, we need to choose $c(x)^3 \notin F$.  So take $c(x) = \tau^{1/9}$.  We then take $f_x(t) = (t-\tau^{1/3})^3$.  Then in $u(L,f)$, we have $(x^3-x-\tau^{1/3})^3=0$, thus $x^9 = x^3+\tau$.  Let $V = u(L,f) \otimes_{u(S,f)} W$.  Then $V$ has basis $\{v_0, \dots, v_8\}$ where $v_r = x^r \otimes w$.  The action of $x$ on $V$ is given by $xv_r = v_{r+1}$ for $r = 0 , \dots, 7$ and $xv_8 = \tau v_0 +v_3$. To calculate the action of $y$, we use the commutation rule $yx^r = (x - 1)^ry$.  We obtain
\begin{align*}
yv_0 &= v_0, & yv_1 &= -v_0 +v_1,\\
yv_2&=v_0+v_1 +v_2, &yv_3&=-v_0+v_3, \\yv_4&=v_0-v_1-v_3+v_4, & yv_5&= -v_0 -v_1-v_2+v_3+v_4+v_5,\\
yv_6&=v_0+v_3+v_6, & yv_7&=-v_0+ v_1 -v_3+v_4 -v_6,\\
yv_8&= v_0+v_1 +v_2 +v_3+v_4 +v_5 + v_6+v_7.\\
\end{align*}
\end{example}

In the above example, $(t-\tau^{1/3})^3$ is the minimal polynomial of the action of $x^3-x^\scp$ on $\bar{V}$.  It follows that the solution spaces in $\bar{V}$ of $(x^3-x-\tau^{1/3})^rv=0$ for $r = 1,2,3$ are distinct.  It follows that $V$ is irreducible and that $\bar{V}$ has a unique composition series with three isomorphic $3$-dimensional factors.

\section{Use of $p$-envelopes.} \label{envelopes}
In this section,  we drop the assumption that $L$ is a restricted Lie algebra with $S$ a $\p$-subalgebra and investigate the extent to which the use of $p$-envelopes can  provide the needed extra structure.  Let $(L^e,\p)$ be a $\p$-envelope of $L$ and let $S_\scp$ be the $\p$-closure of $S$.

Let $\hat{L} \subset U(L)$ be the universal $p$-envelope of $L$ and let $\hat{i} : \hat{L} \to L^e$ be the map given by the universal property of $\hat{L}$.  We choose a cobasis $\{b_1, \dots,b_r\}$ of $L$ in $L^e$.  For each $\alpha$, we choose $a_\alpha \in \hat{i}^{-1}(b_\alpha)$ and put $A = \langle a_1, \dots, a_r, L \rangle$.  Then $\hat{i}|A \to L^e$ is a vector space isomorphism.  But by \cite[Lemma 2.5.5(1)]{SF}, $A$ is a subalgebra of $\hat{L}$.  Since $\hat{i}| A \to L^e$ is the restriction of an algebra homomorphism, it is an algebra isomorphism.  Let $j: L^e \to U(L)$ be the composite of $(\hat{i}|A)^{-1}$ and the inclusion $\hat{L} \to U(L)$.  For any $L$-module $V$, the action extends to an action of $U(L)$ making $V$ a $U(L)$-module.  Restricting this action to $L^e$ makes $V$ into an $L^e$-module.  We thus obtain a functor $J_L$ from $L$-modules to $L^e$-modules.  Similarly, we can construct a functor $J_S$ from $S$-modules to $S_\scp$-modules.

Now suppose that $S_\scp \cap L = S$.  In this case, we can choose the cobasis $\{b_1, \dots,b_r\}$ of $L$ in $L^e$ such that $\{b_1, \dots,b_s\}$ is a cobasis of $S$ in $S_\scp$.  In this case, the two actions are compatible, that is, for $x \in S_\scp$, the action of $x$ on $J_L(V)$ is the same as the action of $x$ on $J_S(\res_S^L(V))$.

Now suppose that we are given a finite set $X$ of $L$-modules.  We construct $\C'$ as in Theorem \ref{Thres}, with $\D'$ the category of all $S_\scp$-modules.  Let $\C$ be the category of $L$-modules $V$ such that $J_L(V) \in \C'$. For an $S$-module $W$, we define $\ind_\D^\C(W)$ to be the $L$-submodule of $\ind_{\D'}^{\C'}(J_S(W))$ generated by $1 \otimes W$.

\begin{lemma} \label{Lemunres} Suppose $(L^e,\p)$ is a \p-envelope of $L$ and that $S_\scp \cap L = S$.  Then $\ind_\D^\C: \D \to \C$ is a left adjoint to $\res_\D^\C :\C \to \D$.
\end{lemma} 

\begin{proof} Let $W$ be an $S$-module and let $V$ be an $L$-module.  Let $f:W \to \res_S^L(V)$ be an $S$-module homomorphism.  Since $J_S$ and $J_L$ are compatible, we have an $S_\scp$-module homomorphism $g: J_S(W) \to \res_{D'}^{\C'}(J_L(V))$.  We then have the corresponding $L^e$-module homomorphism $h: \ind_{\D'}^{\C'}(J_S(W)) \to J_L(V)$, mapping $1 \otimes w$ to $g(w) = f(w) \in \res_\D^\C(V)$.  Since $\ind_\D^\C(W)$ is generated by the $1 \otimes w$ for $w \in W$, the map $h|\ind_\D^\C(W) : \to \res_\D^\C(V)$ is uniquely determined by $f$.
\end{proof}

We now show that for any $L$ and subalgebra $S$, we can choose a $p$-envelope $(L^e, \p)$ of $L$ such that the condition $S_\scp \cap L = S$ is satisfied.  We start with any finite-dimensional $p$-envelope $(L_1,\p_1)$ of $L$.  We choose a basis $\{e_1, \dots,e_n\}$ of $L_1$ such that $\{e_1, \dots, e_s\}$ is a basis of $S_{\scp_1}$. We adjoin an element $z_i$ for $i=1,\dots s$, defining $L^e = \langle z_1, \dots, z_s, L_1\rangle$ with $[z_i,L^e]=0$ and $z_i^{\scp_1}=0$.  We define a new $p$-operation $\p$ on $L^e$ by setting $e_i^\scp = b_i = e_i^{\scp_1} +z_i$ if $i \le  s$ and $e^\scp = e_i^{\scp_1}$ otherwise.

\begin{lemma}\label{psub} With $(L^e,\p)$ constructed as above, $S_\scp \cap L = S$. 
\end{lemma}

\begin{proof} Put $B=\langle b_1, \dots, b_s\rangle$, $Z = \langle z_1, \dots, z_s \rangle$ and $T = S +B$.   Since $b_i = e_i^{\scp_1} +z_i$ and $ e_i^{\scp_1} \in S_{\scp_1}$, any $x \in T$ has the form $x=a+z$ for some $a \in S_{\scp_1}$ and $z \in Z$.  Thus for $x,y \in T$, we have $[x,y] \in (S_{\scp_1})'$.  But $(S_{\scp_1})' = S'$ by \cite[Proposition 2.1.3(2)]{SF}.  Thus $T$ is a subalgebra and $T' = S'$.  Further, for $x = \sum_{i=1}^s(\lambda_i e_i +\mu_i z_i)\in T$, we have $x^\scp = \sum_{i=1}^s \lambda_i^p e_i^\scp +k$ where $k \in T' =  S'$.  It follows that $T$ is a \p-subalgebra of $L^e$, so $T \supseteq S_\scp$.  
Now suppose $x = a + \sum_i \lambda_i b_i \in T$, $a \in S$.  If any $\lambda_i \ne 0$, then $x \notin L$.  Thus $T \cap L = S$.
\end{proof}

\begin{theorem} \label{general} Let $L$ be a Lie algebra over a field $F$ of characteristic $p \ne0$ and let $S$ be a subalgebra of $L$.  Let $X$ be a finite set of finite-dimensional $L$-modules.  Let $D$ be the category of all finite-dimensional $S$-modules.  Then there exists a category $\C$ of finite-dimensional $L$-modules such that $X\subseteq  \C$ and the restriction functor $\res_\D^\C:\C \to \D$ has a left adjoint $\ind_\D^\C:  \to \C$.
\end{theorem}

\begin{proof} By Lemma \ref{psub}, we can choose a $p$-envelope $(L^e,\p)$ of $L$ such that $S_\scp \cap L = S$.  The result follows by Lemma \ref{Lemunres}.
\end{proof}

The following unpleasantly complicated example is designed to illustrate the complications which can arise in dealing with non-restrictable Lie algebras.

\begin{example} Let $L = \langle x, y, d_1, e_1, e_2,e_3, e_4 \rangle$ with the non-zero products $[x, e_1 = e_1+e_2, [x, e_2] = e_2, [d_1, e_1] = e_1, [d_1, e_2]=e_2, [y,e_3]=e_3+e_4$ and $[y,e_4] = e_4$.  Put $S= \langle x, e_1, e_2, e_3, e_4 \rangle$. 
\end{example} 

We take $L_1 = \langle d_2, L \rangle$ with $x^{\scp_1}= d_1, y^{\scp_1} = d_2$ and $d_i^{\scp_1}= e_i^{\scp_i} =0$ for all $i$.  Then $(L_1, \p_1)$ is a minimal \p-envelope of $L$.  But $S_{\scp_1} = \langle d_1, S \rangle$ and $S_{\scp_1} \cap L \ne S$.  We do not need to adjoin central elements corresponding to all elements of the basis of $S_{\scp_1}$ as in Lemma \ref{psub}.  In this example, it suffices to adjoin one central element $z$, forming $L^e = \langle z, L \rangle$.  We set $b = d_1+z$, $x^\scp = b$ and $z^\scp = 0$.  We then have $S_\scp = \langle b, S\rangle$ and $S_\scp \cap L = S$ as required.

We now construct the functor $J$ using the cobasis $\{d_2,b\}$ of $L$ in $L^e$.  We have to choose elements of $\hat{u}(L)$ which map to these cobasis elements under $\hat{i}$.  We can set $j(b) = x^p$.  Since $j(d_1) = d_1$, this gives $j(z) = x^p -d_1$.  Likewise, we can set  $j(d_2) = y^p$.  For any $L$-module $V$ giving the representation $\rho$, we have the action on $JV$ given by $\rho(b) = \rho(x)^p$ and $\rho(d_2)= \rho(y)^p$.  We also have $\rho(z) = \rho(x)^p - \rho(d_1)$.
 
Now suppose that $X$ contains only the one $L$-module $V = \langle v\rangle$ with $xv  =v$ and the other basis elements of $L$ acting trivially.  To construct the $f$-reduced enveloping algebra $u(L,f)$, we must choose polynomials $f_y(t), f_z(t)$ and $ f_{d_2}(t)$  corresponding to the cobasis $\{y, z, d_2\}$ of $S_\scp$ in $L^e$.  We have that  $(y^p-y^\scp)V = 0$, so we can take $f_y(t) = t$.  Likewise, we can take $f_z(t)=f_{d_2}(t)=t$.  In $u(L,f)$, we have $z^p=d_2^p=0$ and $y^p = d_2$.  For any $S$-module $W$, the $L$-submodule of $u(L,f) \otimes_{u(S,f)} W$ generated by $1 \otimes W$ is $\langle y^rd_2^s \otimes W \mid r,s < p\rangle \ne u(L,f) \otimes_{u(S,f)} W$.

\bibliographystyle{amsplain}

\end{document}